\magnification=1100
\baselineskip=14truept
\voffset=.75in  
\hoffset=1truein
\hsize=4.5truein
\vsize=7.75truein
\parindent=.166666in
\pretolerance=500 \tolerance=1000 \brokenpenalty=5000


\def\anote#1#2#3{\smash{\kern#1in{\raise#2in\hbox{#3}}}%
  \nointerlineskip}     
\def\note#1{%
  \hfuzz=50pt%
  \vadjust{%
    \setbox1=\vtop{%
      \hsize 3cm\parindent=0pt\eightrm\baselineskip=9pt%
      \rightskip=4mm plus 4mm\raggedright#1%
      }%
    \hbox{\kern-4cm\smash{\box1}\hfil\par}%
    }%
  \hfuzz=0pt
  }

\font\eightrm=cmr8%
\font\eightbf=cmb8%
\font\eightcmmi=cmmi8%
\font\tenbb=msbm10%
\font\sevenbb=msbm7%
\font\fivebb=msbm5%
\newfam\bbfam%
\textfont\bbfam=\tenbb%
\scriptfont\bbfam=\sevenbb%
\scriptscriptfont\bbfam=\fivebb%
\font\tencmssi=cmssi10%
\font\sevencmssi=cmssi7%
\font\fivecmssi=cmssi5%
\newfam\ssfam%
\textfont\ssfam=\tencmssi%
\scriptfont\ssfam=\sevencmssi%
\scriptscriptfont\ssfam=\fivecmssi%
\def\ssi{\fam\ssfam\tencmssi}%

\font\tenmsam=msam10%
\font\sevenmsam=msam7%
\font\fivemsam=msam5%

\def\bb{\fam\bbfam\tenbb}%

\def\hexdigit#1{\ifnum#1<10 \number#1\else%
  \ifnum#1=10 A\else\ifnum#1=11 B\else\ifnum#1=12 C\else%
  \ifnum#1=13 D\else\ifnum#1=14 E\else\ifnum#1=15 F\fi%
  \fi\fi\fi\fi\fi\fi}
\newfam\msamfam%
\textfont\msamfam=\tenmsam%
\scriptfont\msamfam=\sevenmsam%
\scriptscriptfont\msamfam=\fivemsam%
\mathchardef\leq"3\hexdigit\msamfam 36%
\mathchardef\geq"3\hexdigit\msamfam 3E%

\def\BS{{\ssi S}}
\def\d{\,{\rm d}}
\def\D{{\rm D}}

\long\def\DoNotPrint#1{\relax}

\def\frac#1#2{{#1\over #2}}

\def\Id{{\rm Id}}

\def\L{{\ssi L}}

\def\oF{\overline F{}}
\def\oG{\overline G{}}

\def\oL{\overline L{}}
\def\oW{\overline W{}}
\def\qed{~~~{\vrule height .9ex width .8ex depth -.1ex}}
\def\sign{\hbox{\rm sign}}

\def\II{{\bb I}}
\def\NN{{\bb N}\kern .5pt}
\def\RR{{\bb R}}

%
%
%
%
\def\uncatcodespecials 
    {\def\do##1{\catcode`##1=12}\dospecials}%
{\catcode`\^^I=\active \gdef^^I{\ \ \ \ }
 \catcode`\`=\active\gdef`{\relax\lq}}
\def\setupverbatim 
    {\parindent=0pt\tt %
     \spaceskip=0pt \xspaceskip=0pt 
     \catcode`\^^I=\active %
     \catcode`\`=\active %
     \def\par{\leavevmode\endgraf}
     \obeylines \uncatcodespecials \obeyspaces %
     }%
{\obeyspaces \global\let =\ }
%
%
%
%
%
%

\centerline{\bf TAIL EXPANSIONS FOR THE DISTRIBUTION}
\centerline{\bf OF THE MAXIMUM OF A RANDOM WALK}
\centerline{\bf WITH NEGATIVE DRIFT}
\centerline{\bf AND REGULARLY VARYING INCREMENTS}

\bigskip
 
\centerline{Ph.\ Barbe$^{(1)}$, W.P.\ McCormick$^{(2)}$ and C.\ Zhang$^{(2)}$}
\centerline{${}^{(1)}$CNRS, France, and ${}^{(2)}$University of Georgia}
 
{\narrower
\baselineskip=9pt\parindent=0pt\eightrm

\bigskip

{\eightbf Abstract.}  Let {\eightcmmi F} be a distribution function
with negative mean and regularly varying right tail. Under a mild
smoothness condition we derive higher order asymptotic expansions for
the tail distribution of the maxima of the random walk generated by
{\eightcmmi F}. An application to ruin probabilities is developed.

\bigskip


\bigskip
 
\noindent{\eightbf Keywords:} tail expansion, random walk, regularly
varying, Wiener-Hopf factor, ruin probability

}

\bigskip\bigskip


\noindent{\bf 1.\ Introduction.\ }
There is hardly a more basic stochastic model than a random walk, and
for random walks with negative drift, a basic issue of study is the
distribution of its global maximum. One reason for interest in this
quantity is its connection to queueing processes. For a GI/G/1 queue,
which is stable in the sense that the mean interarrival time exceeds
the mean service time, the waiting time that the $n$-th arriving
customer needs to wait until service begins has a limiting
distribution as $n$ tends to infinity. This is given by the
distribution of the global maximum of a random walk with negative
drift; see Asmussen (1987, \S III.7). In an insurance-risk setting, the
distribution of the global maximum of a random walk with negative
drift directly appears in computing ruin probabilities over an
infinite horizon; see, for example, Embrechts, Kl\"uppelberg and Mikosch
(1997, \S 1.1).

Let ${(X_n)}_{n \ge 1}$ be a sequence of independent and identically
distributed random variables having negative mean. The associated
random walk is defined by $S_0=0$ and for any integer $n$ positive,
$S_n=X_1+\cdots +X_n$. The distribution of its maximum, $M=\max_{n\geq
0}S_n$, can be represented as a compound-geometric
distribution as follows. We first agree that the minimum of the
empty set is $+\infty$. Then, let $\tau$ denote the hitting time for the 
positive half-line
$$
  \tau = \min\{\,n:S_1 \le 0, \dots ,S_{n-1} \le 0\, ,\, S_n >0 \,\} \, .
$$
This hitting time may be infinite, but it is finite with probability
$$
  p
  =P\{\,\tau <\infty\,\} 
  = 1-P\{\, S_1\le 0, S_2 \le 0, \dots \,\} \, .
$$
Recall that the first strict ascending ladder height distribution is
defined by 
$$
  F_+(x) = P\{\,S_\tau \le x\, ,\,\tau < \infty\,\} \, .
$$ 
Since the random walk has a negative drift,
$F_+$ is a defective distribution with defect $1-p$. By the Sparre Andersen 
identity (Feller, 1971, \S XII.7) and Abel's lemma (Karlin, 1975, \S II.5),
$$
  p = 1 -\exp\Bigl(-\sum_{n\geq 1} {1\over n} P\{\,S_n >0\,\}\Bigr) \, ; 
$$
see also Chung (1974). It follows that $M$ has a compound-geometric
distribution subordinate to the distribution $H=p^{-1}F_+$ and with
subordinator a geometric distribution with parameter $p$. More explicitly
and following Feller (1971, \S XII.5), writing $H^{\star n}$ the $n$-fold 
convolution of $H$, the distribution $W$ of $M$ is
$$
  W=(1-p)\sum_{n\geq 0}p^nH^{\star n}\, .
  \eqno{(1.1)}
$$

This step has replaced the original question of analyzing the
distribution of global maximum of a random walk with the more
elementary question of analyzing that of a compound sum, at the price,
however, of introducing a derived distribution, namely, the ascending
ladder height distribution, which requires its own analysis. In the
case of a heavy-tailed step size distribution as prescribed by a
subexponentiality assumption, Veraverbeke (1977) supplies an answer to
this question through use of the distributional form of the
Wiener-Hopf factorization. His result establishes inheritability of
the subexponential property of the right Wiener-Hopf factor, i.e.\ that
factor having mass concentrated in positive half line, from that of
the underlying distribution. To state that result, we agree that for
any possibly defective distribution function $G$, we write $\oG$ its
tail, that is the function whose value at $x$ is
$$
  \oG(x)=\lim_{t\to\infty} G(t)-G(x) \, .
$$
Consider the Wiener-Hopf factorization $F=F_++F_--F_+\star F_-$, where
$F_-$ and $F_+$ are concentrated on $(-\infty,0\,]$ and $(0,\infty)$
respectively (see Feller, 1971, \S XII.3). Let $\mu$ be the mean of $F$, 
which we assume to be negative. Veraverbeke (1977) shows that, as $x$ 
tends to infinity,
$$
  \oF_+(x)
  \sim {1-p\over -\mu}\int_x^\infty \oF(t)\d t \, .
$$ 
We remark that the Wiener-Hopf factors are given by the strict
ascending ladder height distribution and the weak descending ladder
height distribution.

With these two steps in place, a first-order analysis of the
distribution of $M$ may be completed by using a result on
tail-area asymptotics for subordinated probability distributions --- in
this case for the subordinator given by a geometric distribution. For
example, the result in Athreya and Ney (1972, \S IV.4) for first-order
asymptotics of compound subexponential distributions with geometric
subordinator gives the expected result that, under the
assumption of $F$ subexponential with negative mean $\mu$,
$$ 
  \oW(x) 
  \sim {1\over 1-p}\oF_+(x) 
  \sim {-1\over\mu}\int_x^\infty \oF(t)\d t\, ,
$$
as $x$ tends to infinity.
This last step of analysing tail areas for subordinated
distributions is a subject of much interest in the literature. In the
case of subexponential subordinate distributions, we mention the
paper, Embrechts, Goldie and Veraverbeke (1979), 
who prove that for any
integer valued random variable $N$ independent of the 
sequence ${(X_i)}_{i\geq 1}$ and such that $Ez^N$ is analytic at $z=1$,
$$
  P\{\,S_N >x\,\} \sim EN\, P\{X_1 >x\}
$$ 
as $x$ tends to infinity.

Finally, we mention that when $\oF$ is regularly varying with index
in the range from $-1$ to $-2$ and the mean is finite and negative,
Omey and Willekens (1986) establish a second-order result for the tail
$\oW$; see also Geluk (1992, 1996). 

In a different range of tail heaviness which we will not explore
in this paper, Feller (1971, \S XII.5, Example c) considers the tail behavior
of $W$ when $F$ has a moment generating function finite in a neighborhood
of the origin. He shows that if the moment generating function of $F$
is $1$ at some positive $\kappa$, and if the number $\beta 
= \int_0^\infty x \exp(\kappa x) \d F_+(x)$ is finite, then
$$
  \oW(x) \sim {1-p\over \beta \kappa} e^{-\kappa x}
$$
as $x$ tends to infinity. If $\beta$ is infinite, the result should be read as
$\oW =o(e^{-\kappa x})$.

\bigskip


\noindent{\bf 2.\ Expansions.\ }
In this paper, we are interested in proving higher-order tail area
asymptotics for the distribution $W$ of the maximum $M$ of the random
walk. As previously noted,
the representation of this distribution as a compound-geometric allows
a restatement of this problem as one of establishing higher-order
results for certain compound-geometric distributions. We remark that
Omey and Willekens (1987) provide second-order results for such a
compound distribution subordinate to a distribution with finite first
moment which satisfies certain smoothness and regularity conditions,
e.g. membership in a subclass of subexponential distributions
including distributions with regularly varying tails with index of
regular variation at most $-1$. Note however that their results cannot
be automatically applied to derive second-order behavior for
$M$. This is because their result requires certain smoothness
conditions which would need to be established for the right
Wiener-Hopf factor under appropriate conditions on the underlying
distribution.

Our approach to this question is to invoke a result from Barbe and
McCormick (2004) on asymptotic expansions for tail areas of
compound sums. Their Theorem 4.4.1 provides an $m$-term expansion for
tail area, where the number of terms allowed is constrained by
smoothness and moment conditions on the underlying distribution. To
that end, we now present a smoothness condition needed for that
result.

\bigskip

\noindent{\bf Definition.\ }{\it
  A real measurable function $f$ is smoothly varying with index
  $-\alpha$ and order $m$ if it is ultimately $m$-times continuously
  differentiable and the $m$-th derivative $f^{(m)}$ is regular varying
  with index $-\alpha-m$. We denote the set of all such functions by 
  $SR_{-\alpha,m}$.
}

\bigskip

All distributions with regularly varying tails used in applications
are smoothly varying of arbitrary order. Examples include the Pareto,
Cauchy, Student, Burr and log-gamma distributions. Any function smoothly
varying in the sense of Bingham, Goldie and Teugels (1984, \S 1.8.1)
is smoothly varying of any fixed order.

The class $SR_{-\alpha,m}$ may be extended to noninteger orders.  This
is useful to present sharp results. To define $SR_{-\alpha,\omega}$
where $\omega$ is a positive real number, we introduce the following
notation. For any function $h$, let
$$
  \Delta^r_{t,x}(h) 
  =\sign(x) {h\bigl(t(1-x)\bigr) - h(t)\over |x|^r \, h(t)} \, .
$$

\bigskip

\noindent{\bf Definition.} {\quad\it
  Let $\omega$ be a positive real number. Write $\omega = m+r$ where $m$
  is the integer part of $\omega$ and $r$ is in $[\, 0,1)$. A 
  function $h$ is smoothly varying of
  index $-\alpha$ and order $\omega$ if it belongs to $SR_{-\alpha,m}$ and
  $$ 
    \lim_{\delta \rightarrow 0} \limsup_{t \rightarrow \infty}
    \sup_{0<|x| \le \delta}\Delta^r_{t,x}(h)=0 \, .
  $$ We write $SR_{-\alpha,\omega}$ for the class of all such functions.  
}

\bigskip  

We remark that the spaces $SR_{-\alpha,\omega}$ are nested, for $SR_{-\alpha,r}
\supset SR_{-\alpha,s}$ for $r < s$. In particular, if $\omega$ is positive
with integer part $m$ and $r=\omega-m$, membership in $SR_{-\alpha,\omega}$ is
guaranteed by that in $SR_{-\alpha,m+1}$, that is by checking that the
$m+1$-derivative is regularly varying of index $-\alpha-m-1$.  For further
properties of smoothly varying functions of finite order, we refer to
Barbe and McCormick (2004).

\bigskip

We now introduce the algebraic formalism to express the result. To
that end let $\RR_m[\,\D\,]$ denote the ring of real polynomials in
$\D$ modulo the ideal generated by $\D^{m+1}$. In other words, any
polynomial in $\D$ divisible by $\D^{m+1}$ is set equal to $0$.

For a possibly defective distribution function $G$ with at least $k$ 
moments finite, we write $\mu_{G,k}$ its $k-$th moment. Note in
particular that $\mu_{G,0}$ is the total mass of $G$, equal to $1$ if and
only if $G$ is not defective.

\bigskip

\noindent{\bf Definition.\ }{\it
  The Laplace character of order $m$ of a possibly defective distribution $G$
  having a finite $m$-th moment is the element of $\RR_m[\,\D\,]$ given by
  $$ 
    \L_{G,m} = \sum_{0 \le k \le m}{(-1)^k\over k!} \mu_{G,k}\D^k \, .
  $$ 
}

Note that since the map which associates to a measure its $k$-th moment
is linear on its domain, the map $G\mapsto \L_{G,m}$ is linear on its domain.

The backward signed shift $\BS$ on polymomials in $\D$ is defined linearly by
$\BS\D^0=0$ and whenever $j$ is a positive integer, $\BS\D^j=-\D^{j-1}$.
It maps $\RR_m[\,\D\,]$ to $\RR_{m-1}[\,\D\,]$.

We define the inverse of the differentiation on some functions as follows. 
If $f$ is a function regularly varying of index less than $-1$, we set
$$
  \D^{-1}f(t)=-\int_t^\infty f(x)\d x \, .
$$
Clearly, $\D\D^{-1}$ is the identity on functions which are regularly
varying of index less than $-1$, while $\D^{-1}\D$ is the identity on 
the smoothly varying functions of negative index and order at least $1$.

We now present our main result. Recall $F$ denotes the step size
distribution about which we assume its first moment is negative and that
its right tail is regularly varying of index $-\alpha$.
The strict ascending ladder height distribution is $F_+$ and
$W$ is the distribution of $M$. Let 
$$
  \kappa=\sup\Bigl\{\, r\geq 0 \, :\, \int_{-\infty}^0 |x|^r \d F(x) <\infty 
  \,\Bigr\} \, .
$$
This may be less than $\alpha$ if the lower tail
of $F$ is heavier than the upper one. In the following theorem, it is
implicitly supposed that $\kappa$ is greater than $1$.

\bigskip

\noindent{\bf Theorem.\ } {\it
  Suppose that $\oF$ is smoothly varying of index $-\alpha$ and order
  $\omega$. Then, for any integer $m$ at least $1$ and less than 
  $\omega\wedge\alpha\wedge \kappa$, the moments $\mu_{F_-,m}$ and
  $\mu_{F_+,m-1}$ are finite
  and
  $$
    \oW = (1-p)
    (\Id - \L_{F_+,m-1})^{-2}(\BS\L_{F_-,m})^{-1}(\D^{-1}\oF) 
    +o(\Id^{-m+2} \oF)\, .
  $$ 
}

\noindent {\bf Remark.\ } As previously mentionned, Laplace characters
of order $m-1$ are elements of the ring $\RR_{m-1}[\,\D\,]$. The
inverses $(\Id-L_{F_+,m-1})^{-2}$ and $(\BS \L_{F_-,m})^{-1}$ are taken in 
that ring, multiplied together in that ring, and applied to $\D^{-1}\oF$.

\bigskip

\noindent{\bf Remark.\ } The result may seem a little mysterious and not
so explicit at a first glance. However, the computations related to
Laplace characters can be implemented with a computer algebra package. For
instance, the following very short {\tt Maple} code calculates
the expansion given in the Theorem. In that code, {\tt Fp} and {\tt Fm}
stand for $F_+$ and $F_-$.

\verbatim@
restart; m:=4: mu[Fp,0]:=1-q:
LFp:=sum('(-1)^j*mu[Fp,j]*x^j/j!','j'=0..m-1):
SLFm:=sum('(-1)^j*mu[Fm,j+1]*x^j/(j+1)!','j'=0..m-1):
a:=taylor((1-LFp)^(-2),x=0,m-1):
b:=taylor(SLFm^(-1),x=0,m-1):
expand(convert(q*taylor(a*b,x=0,m-1),polynom)/x);

@

\noindent
One simply replaces {\tt x} and {\tt q} in the output by $\D$ and $1-p$, with 
the convention that $1/x$ should be replaced by $\D^{-1}$.
For instance, using that $\mu_{F,1}=(1-p)\mu_{F_-,1}$, taking $m$ to 
be $4$, we deduce the $3$-terms expansion
$$\eqalign{
  \oW={}
  & {1\over\mu_{F,1}}\D^{-1}\oF 
    + {1\over 2\mu_{F,1}^2} 
    \bigl((1-p)\mu_{F_-,2}-4\mu_{F_+,1}\mu_{F_-,1}\bigr) \oF \cr
  &{} + {1\over 12\mu_{F,1}^3} \Bigl(
    3(1-p)^2\mu_{F_-,2}^2
    +12\mu_{F,1}(\mu_{F_-,1}\mu_{F_+,2}-\mu_{F_+,1}\mu_{F_-,2})\cr
  &\hskip 60pt {}+ 36\mu_{F_-,1}^2\mu_{F_+,1}^2  - 2(1-p)\mu_{F,1}\mu_{F_-,3}
   \Bigr)\oF' \cr
  &{}+o(\Id^{-1}\oF) \, . \cr}
$$

\noindent{\bf Remark.\ } An important point to mention with regard to the
main result is that the expansion it provides for the tail distribution
$\oW$ is based on the underlying distribution $F$ of the random walk, its
derivatives and its integrated tail. This is notable, because the starting
point to obtain this result is that of a tail area expansion for a subordinated
distribution based on underlying distribution given by $F_+$. Since $F_+$ is
generally unattainable in an explicit form, the formulation of our main
result is more attractive than that which results from a direct application
of a tail area result for subordinated distributions. It is a comment on
the usefulness of this algebraic approach that such an improvement is
so easily and transparently attained compared to the effort to accomplish
the same goal analytically. We conclude this remark by noting that our
proof shows that a penultimate expansion based on $F_+$ is given by
$$
  \oW=(1-p)(\Id-\L_{F_+,m})^{-2}\oF_+ + o(\Id^{-m+1}\oF)
  \eqno{(2.1)}
$$
provided $m$ is less than $\omega\wedge (\alpha-1)$.
When $m$ is $1$ in the above, we obtain a second-order result in agreement with
Theorem 2.2 in Omey and Willekens (1987). Comparing this formula with that
given in our theorem, we see that the latter is slightly less accurate.
The reason is that, under the assumption of the theorem, replacement of
$\oF_+$ with an approximation based on $\oF$, its derivatives and its integral
comes with a one-order lower error bound, viz.
$$
  \oF_+=(\BS L_{F_-,m})^{-1}\D^{-1}\oF + o(\Id^{-m+2}\oF) \, .
$$

\bigskip

\noindent{\bf Application.\ }
Finally, we present an application to insurance risk. To that end, we
introduce some notation. Let $R_0=x$ be the initial capital of an
insurance company. We assume that the claim amounts, ${(A_n)}_{n\geq 1}$, 
are independent, with common distribution function $L$ having
a smoothly varying tail of index $-\alpha$ and order $m+1$. We also
assume that the interclaim
times ${(T_n)}_{n\geq 1}$ are independent, with common distribution $K$, and 
independent of the claim amounts. Finally, we assume the intensity
of the gross risk premium is some positive $c$. The net
loss to the company in period $n$ is $X_n=A_n-cT_n$. The sequence 
${(X_n)}_{n \geq 1}$ is a sequence of independent random variables with 
distribution
$F(x)=\int_{0}^{\infty}L(x+ct)\d K(t)$. Under the assumptions
on $L$, it follows that $F$ is ultimately $m$-times differentiable and
$$
  {F^{(m)}(x) \over L^{(m)}(x)}
  =\int_0^\infty\frac{L^{(m)}(x+ct)}{L^{(m)}(x)}\d K(t)\, .
$$ 
This implies that $\oF$ is smoothly varying of index $-\alpha$ and order $m$.
Let $\psi (R_0)$ be the probability of eventual ruin given $R_0$.
Writing as before $S_n$ for the random walk with increment $X_i$ and $M$
for its maximum, $\psi(x) = P\{\,M> x\,\}$. We follow our
established notation and set $F_+$ and $F_-$ for the strict ascending and
weak descending ladder height distributions for the
random walk $S_n$. We assume that $X_1$ has negative expectation.  
We have following expansion of the ruin probability, which obviously
follow from the theorem.

\bigskip

\noindent{\bf Corollary.\ }{\it
  Assume that $\oL$ is smoothly varying of index $-\alpha$ and order $m+1$.
  Assume also that
  $\mu=EA_1-cET_1$ is finite and negative and that $m$ less than
  $\alpha$. Then,
  $$\displaylines{\qquad
    \psi(x)
    =(1-p) (\Id - \L_{F_+,m-1})^{-2}(\BS\L_{F_-,m})^{-1}
    \D^{-1}\oF(x) 
    \hfill\cr\noalign{\vskip 2pt}\hfill
    {}+o\bigl(x^{-m+2}\oF(x)\bigr) \, .\qquad\cr}
  $$ 
}



\noindent{\bf 3. Proof of the theorem.\ } We first show that $\mu_{F_-,m}$
is indeed finite. The distributional form of the Wiener-Hopf factorization
implies that on the negative half-line,
$$
  F_-=F+F_+\star F_- \, .
$$
Since $F_+$ has defect $1-p$, this yields $F_-\leq F+pF_-$, from which
we deduce $F_-\leq (1-p)^{-1}F$. This proves the finiteness of $\mu_{F_-,m}$.

We begin the proof of the main part of our theorem by establishing a 
preparatory lemma. 

\bigskip

\noindent{\bf Lemma 1.\ } {\it 
   Let ${(Y_i)}_{i\geq 1}$ be a sequence of nonnegative
   random variables, independent and identically distributed with finite
   and positive mean. Let
   ${(Z_n)}_{n\geq 0}$ be their corresponding random walk. 
   Furthermore, let $f$ be a regularly varying function of index $-\beta$ less 
   than $-1$. Then,
   $$
     \lim_{x\to\infty} {1\over xf(x)}
     \sum_{n\geq 0} E f(x+Z_n) = {1\over (\beta-1)EY_1} \, .
   $$ 
}

\noindent
{\bf Proof.\ } We may assume that $f$ is ultimately positive.  Since
$f$ is regularly varying of negative index, it is asymptotically
equivalent to a nonincreasing function (Bingham, Goldie and Teugels,
\S 1.5.2).  Therefore, we can assume without any loss of generality
that $f$ is nonincreasing. Let $\theta$ be positive and less than the
mean of the $Y_i$'s. We have the trivial bound
$$
  \sum_{n\geq 0} Ef(x+Z_n)
  \leq \sum_{n\geq 0}f(x+\theta n ) 
  + f(x)\sum_{n\geq 0}P\{\, Z_n\leq \theta n\,\}
  \, .
$$
Since $f$ is regularly varying,
$$
  \sum_{n\geq 0} f(x+\theta n)\sim \int_0^\infty f(x+\theta s)\d s \, .
$$
as $x$ tends to infinity. The change of variable $\theta s=xz$ yields
$$\eqalign{
  \int_0^\infty f(x+\theta s)\d s
  &{}= {x\over \theta} \int_0^\infty f\bigl( x(1+z)\bigr) \d z \cr
  &{}\sim {xf(x)\over \theta} \int_0^\infty (1+z)^{-\beta}\d z \cr
  &{}={xf(x)\over \theta (\beta-1)} \, . \cr
  }
$$
Let $M$ be such that $E(Y_1\wedge M)$ is greater than $\theta$. Such
$M$ exists by monotone convergence of $Y_1\wedge M$ to $Y_1$. Let
$(Z_n^M)_{n\geq 0}$ be the random walk associated to the sequence
${(Y_i\wedge M)}_{i\geq 1}$.  By the Hsu and Robbins theorem (see Chow
and Teicher, 1988, \S 10.4), the series $\sum_{n\geq 0}P\{\, Z_n^M\leq
n\theta\,\}$ is finite. Since $Z_n^M$ is at most $Z_n$, this series is
at least $\sum_{n\geq 0}P\{\, Z_n\leq n\theta\,\}$, and the latter is
finite as well.  Therefore, since $\theta$ is any positive number less
than $EY_1$,
$$
  \limsup_{x\to\infty} {1\over xf(x)} \sum_{n\geq 0} Ef(x+Z_n)
  \leq {1\over (\beta-1) EY_1}\, .
$$

To obtain a matching lower bound, let now $\theta$ be a number greater 
than $1$, and let $\epsilon$ be a positive real number. Since $f$ is
ultimately positive and nonincreasing, $\sum_{n\geq 0}Ef(x+Z_n)$ is
ultimately at least
$$
  \sum_{n\geq 0} Ef\bigr(x+n(EY_1+\epsilon)\bigl)
  \II\{\, |Z_n-nEY_1|\leq \epsilon n\, ;\,
              \theta^{-1}x\leq n EY_1\leq \theta x\,\} \, .
$$
But if $\theta^{-1}x\leq nEY_1\leq \theta x$, as $x$ tends to infinity,
$$
  f\bigl(x +n(EY_1+\epsilon)\bigr)
  \sim f(x) \Bigl( 1+{n\over x}(EY_1+\epsilon)\Bigr)^{-\beta} \, .
$$
Moreover, in that range of $n$, for $x$ large enough, the strong law
of large numbers implies that $P\{\,|Z_n-nEY_1|\leq n\epsilon\,\}
\geq 1-\epsilon$. Therefore, $\sum_{n\geq 0}Ef(x+Z_n)$ is ultimately at least
$$\displaylines{
  (1-\epsilon) f(x)\sum_{n\geq 0} 
  \Bigl( 1+{n\over x}(EY_1+\epsilon)\Bigr)^{-\beta}
  \II\{\, \theta^{-1}x\leq nEY_1\leq \theta x\,\}
  \hfill\cr\hfill
  \eqalign{
    {}\sim{}& (1-\epsilon)f(x)\int_{\theta^{-1}x/EY_1}^{\theta x/EY_1}
              \Bigl( 1+{s\over x}(EY_1+\epsilon)\Bigr)^{-\beta} \d s \cr
    {}={}   & (1-\epsilon) f(x){1\over 1-\beta} \, {x\over EY_1+\epsilon}
              \biggl[\, \Bigl( 1+{s\over x} (EY_1+\epsilon)\Bigr)^{-\beta+1}
              \,\biggr]_{\theta^{-1}x/EY_1}^{\theta x/EY_1} \, .\cr}
  \cr}
$$
Since $\theta$ and $\epsilon$ are arbitrary, we can make $\theta$ tend to 
infinity after taking the asymptotic equivalent of the lower bound as
$x$ tends to infinity, proving that
$$
  \sum_{n\geq 0} Ef(x+Z_n)\geq \bigl( 1+o(1)\bigr) {xf(x)\over EY_1 (\beta-1)}
$$
as $x$ tends to infinity.\hfill\qed

\bigskip

Note that with unimportant and additional assumptions an alternate proof
of Lemma 1 based on the renewal theorem may be given, slightly shorter, but
not as direct. To sketch it, write $G$ the distribution function of $Y_i$ and
consider the renewal function $U=\sum_{n\geq 0} G^{\star n}$. We see that
$\sum_{n\geq 0} Ef(x+Z_n)=\int f\d U$. When $f$ is smooth, an integration 
by parts and a change of variable bring this integral to the form
$s\int_1^\infty f'(xs)U\bigl( x(s-1)\bigr) \d s$. The renewal theorem
(Feller, 1971, \S XI.3) yields $U\bigl( x(s-1)\bigr)\sim x(s-1)/EY_1$ as
$x$ tends to infinity, uniformly in $s$ at least $1$. The result then follows
by standard arguments involving regular variation.

\medskip

The main argument for proving our theorem is to show
that $\oF_+$ is smoothly varying of index $-\alpha+1$ and same order as $\oF$.
This is stated in the next lemma.

\bigskip

\noindent{\bf Lemma 2.\ } {\it 
  The strict ascending ladder height distribution $F_+$ is smoothly
  varying of index $-\alpha+1$ and same order $\omega$ as $F$.
  Moreover, 
  $$
    \oF_+\sim {-1\over (\alpha-1)\mu_{F_-,1}}\,\Id\, \oF \, .
  $$
}

\noindent{\bf Proof.\ } The proof has four steps.

\noindent {\it Step 1. A representation for $\oF_+$.}
By the distributional form of Wiener-Hopf factorization, we have
$$
  F=F_+ + F_- -F_+ \star F_- \, .
  \eqno{(3.1)}
$$
It is convenient to introduce the following integral operator,
$$
  U_{F_-}g(t)=\int_{-\infty}^0 g(t-u) \d F_-(u) \, .
$$  
As usual, powers of operators are defined inductively. In particular,
$U_{F_-}^0$ is the identity and $U_{F_-}^{n}=U_{F_-}\circ
U_{F_-}^{n-1}$ for any integer $n$ positive.  On $(0,\infty)$, we can 
write (3.1) as $\oF_+ = \oF +\overline{F_+ \star F_-}$, which leads to
$$
  \oF_+= \oF + U_{F_-} \oF_+ \, .
  \eqno{(3.2)}
$$ 
By recursion this yields
$$ 
  \oF_+=\sum_{0\leq i\leq n} U_{F_-}^i \oF + U_{F_-}^{n+1}\oF_+\, .
$$
Note that $F_-$ cannot be the distribution degenerate at $0$ since $F$ 
is assumed to have a negative mean.
Let ${(Y_i)}_{i\geq 1}$ be a sequence of independent random variables,
all with the same distribution $F_-$, and let ${(Z_n)}_{n\geq 0}$ be
their random walk (note that the signs are changed compared to the
previous lemma).  Observe that $\sum_{0\leq i\leq n} U_{F_-}^i \oF$ is
nondecreasing in $n$ and that, by dominated convergence, 
$U_{F_-}^{n+1} \oF_+(x)= E\oF_+(x-Z_{n+1})$
tends to $0$ as $n$ goes to infinity. Consequently, we obtain the
representation
$$
  \oF_+=\sum_{i\geq 0} U_{F_-}^i \oF\, .
  \eqno{(3.3)}
$$
Note that combined with Lemma 1, this representation yields Veraverbeke's
(1977) theorem asserting that
$$
  \oF_+\sim {-1\over (\alpha-1)\mu_{F_-,1}}\,\Id\,\oF \, .
  \eqno{(3.4)}
$$

\noindent{\it Step 2. A representation for $\oF_+^{(k)}$.}
Let $k$ be a positive integer at most $\omega\wedge (\alpha-1)$. 
Using the mean value theorem, there
exists a sequence of real numbers, ${(\theta_n)}_{n\geq 0}$, nonnegative and
at most $1$, such that
$$
\displaylines{
   \sum_{n\geq 0} \Bigl| {1\over\epsilon}\Bigl( U_{F_-}^n\oF^{(k-1)}
   (x+\epsilon)-U_{F_-}^n \oF^{(k-1)}(x)\Bigr)
   -U_{F_-}^n \oF^{(k)}(x)\Bigr|
  \hfill\cr\hfill
  \eqalign{
  {}={}   &\sum_{n\geq 0} \bigl| U_{F_-}^n\oF^{(k)}(x+\theta_n\epsilon)
           -U_{F_-}^n\oF^{(k)}(x)\bigr|\cr
  {}\leq{}&\sum_{n\geq 0} E\bigl|\oF^{(k)}(x+\theta_n\epsilon-Z_n)
                                 -\oF^{(k)}(x-Z_n)\bigr|\, .
        \cr}
  \cr}
$$
Since the absolute value of a difference is at most the sum of the
absolute values, Lemma 1 shows that the above series is bounded
as a function of $x$ and uniformly in $\epsilon$ in some interval
$(0,\eta)$. Moreover, every summand tends to $0$ as $\epsilon$
tends to $0$. Therefore, the series tends to $0$ as $\epsilon$ tends
to infinity. This proves that
$$
  \oF_+^{(k)}=\sum_{i\geq 0} U_{F_-}^i\oF^{(k)}
  \eqno{(3.5)}
$$
on some neighborhood of infinity.

\smallskip

\noindent{\it Step 3.  $\oF_+^{(k)}$ is regularly varying.}  The
asymptotic equivalence in (3.4) implies that $\oF_+$ is regularly
varying with index $-\alpha+1$.  Recall that $k$ is at most
$\omega\wedge (\alpha-1)$.  By assumption $\oF^{(k)}$ is regularly
varying. By representation (3.5) and Lemma 1, $\oF_+^{(k)}$ is regularly
varying of index $-\alpha-k+1$. Taking $k$ to be $m$, that is $\lfloor
\omega\rfloor$, this proves that $\oF_+$ is smoothly varying of index
$-\alpha+1$ and order $m$.

\smallskip

\noindent{\it Step 4. Concluding the proof of the lemma.} 
Following Barbe and McCormick (2004), for a function $h$ define
$$
  \overline\Delta_{\tau,\delta}^r(h)=
  \sup_{t\geq \tau}\sup_{0<|x|\leq\delta} |\Delta_{t,x}^rh| \, .
$$
This quantity is nonincreasing in $\tau$ and nondecreasing in $\delta$. 
Using representation (3.5), we see that
$$\displaylines{\qquad
  \oF_+^{(m)}\bigl( t(1-x)\bigr) -\oF_+^{(m)}\bigl(t\bigr)
  \hfill\cr\hfill
  {}=\sum_{n\geq 0} E\Bigl( F^{(m)}\bigl( t(1-x)-Z_n\bigr) 
  -F^{(m)}(t-Z_n)\Bigr)\, .\qquad\cr}
$$
Consider $x$ in the range $[\, -\delta,\delta\,]\setminus \{\, 0\,\}$. 
Factoring $t-Z_n$ in $t(1-x)-Z_n$, the $n$-th summand in the series above is
at most
$$
  \Bigl({t\over t-Z_n}|x|\Bigr)^r |\oF^{(m)}(t-Z_n)| \;
  \overline\Delta_{t-Z_n,t\delta/(t-Z_n)}^r\oF^{(m)} \, .
$$
Consequently, for $|x|$ positive and at most $\delta$,
$$
  |\Delta_{t,x}^r\oF_+^{(m)}|
  \leq E\sum_{n\geq 0}\Bigl| {\oF^{(m)}(t-Z_n)\over \oF_+^{(m)}(t)}\Bigr|
  \;\overline\Delta_{t,\delta}^r\oF^{(m)} \, .
$$
It follows from step 3, $\oF_+^{(m)}\asymp \Id \oF^{(m)}$. It then
follows from Lemma 1 and our assumption on $F$ that
$$
  \lim_{\delta\to 0}\lim_{t\to\infty}\sup_{0<|x|<\delta} 
  |\Delta_{t,x}^r\oF_+^{(m)}| = 0 \, ,
$$
proving the smooth variation of order $\omega$ of $\oF_+$.\hfill\qed

\bigskip

Finally, we present a technical lemma of some independent interest, 
particularly in the light of Mari\'c's (2000) work. It is needed for the
proof of our main result. We remark that the result is not proved under 
optimal conditions.

\bigskip

\noindent{\bf Lemma 3.\ } {\it
  Let $(a_i)_{0\leq i\leq m}$ be a sequence of real numbers with $a_0$ 
  different from $0$. For any nonnegative integer $k$ at most $m$, define
  the differential operators $P_k(\D)=\sum_{0\leq i\leq k}a_i\D^i$.
  Let $\psi$ be a function.
  Let $f$ and $g$ be two functions smoothly varying with index
  $-\alpha$ and order at least $m$ satisfying the differential equations
  $$
    P_{m-k}(\D)\D^kf=\D^kg +o(\psi)\, 
    ,\qquad \hbox{ $k=0,1,\ldots ,m$.}
  $$
  Then, viewing $P_m(\D)$ in $\RR_m[\,\D\,]$,}
  $$
    f=P_m(\D)^{-1}g+o(\psi) \, .
  $$

\noindent
The lemma may be interpreted as saying that if the
functions $\D^kg$ have a generalized asymptotic expansion in the the
asymptotic scale $\D^kf$, then $f$ has a generalized asymptotic
expansion in the asymptotic scale $\D^kg$.

\bigskip

\noindent{\bf Proof.} Write $b_k$ the $k$-th coefficient of $P_m(\D)^{-1}$.
Then 
$$
  P_m(\D)^{-1}g
  =\sum_{0\leq k\leq m}b_k \D^kg \, .
  \eqno{(3.6)}
$$
In this sum, by assumption, we can replace $\D^kg$ by
$P_{m-k}(\D)\D^kf+o(\psi)$. Since $P_{m-k}(\D)\D^k=P_m(\D)\D^k$ in 
$\RR_m[\,\D\,]$, the definition of the $b_k$ and (3.6) yield
$P_m(\D)^{-1}g=f+o(\psi)$, which is the result.\hfill$\qed$

\bigskip

We now conclude the proof of the main theorem.
Using (3.2) and applying Lemma 2 and a variant of Theorem 2.3.1 in Barbe and 
McCormick (2004), we 
obtain for any nonnegative $k$ at most $m$,
$$\eqalign{
  \oF^{(k)}
  &{}=\oF_+^{(k)} -\L_{F_-,m-k}\oF_+^{(k)} +o(\Id^{-m}\oF_+)\cr
  &{}=\BS\L_{F_-,m-k}\D^{k+1}\oF_+ + o(\Id^{-m}\oF_+)\, .\cr}
$$
By Veraverbeke's (1977) theorem or (3.4), this implies
$$
  \BS\L_{F_-,m-k}\D^k\D\oF_+ 
  = \D^k\oF+o(\Id^{-m+1}\oF) \, .
$$
Applying Lemma 3, we obtain
$$
  \D\oF_+= (\BS\L_{F_-,m})^{-1}\oF +o(\Id^{-m+1}\oF) \, .
$$
Hence, integrating,
$$
  \oF_+=(\BS\L_{F_-,m})^{-1}\D^{-1}\oF +o(\Id^{-m+2}\oF) \, .
  \eqno{(3.7)}
$$
Representation (1.1) and Theorem 4.4.1 in Barbe and McCormick (2004),
upon noting that if $N$ is a random variable with geometric
distribution with parameter $p$, then 
$EN\L_{H,m}^{N-1}$ is $p(1-p)(\Id-p\L_{H,m})^{-2}$ in $\RR_m[\,\D\,]$,
yield formula (2.1), that is
$$
  \oW=(1-p)(\Id-\L_{F_+,m})^{-2}\oF_+ +o(\Id^{-m}\oF_+)\, .
$$
To obtain the statement of the Theorem, we again use representation (1.1) 
and apply Theorem 4.4.1 in Barbe and McCormick (2004) to obtain that if
$m$ is less than $\alpha\wedge\omega$,
$$
  \oW=(1-p)(\Id-\L_{F_+,m-1})^{-2}\oF_+ + o(\Id^{-m+1}\oF_+) \, .
$$
Then, we use Lemma 2 and (3.7) to conclude.\hfill\qed

\bigskip


\noindent{\bf References}
\medskip

{\leftskip=\parindent \parindent=-\parindent
 \par

K.B.\ Athreya, P.\ Ney (1972). {\sl Branching Processes}, Springer.

S.\ Asmussen (1987). {\sl Applied Probability and Queues}, Wiley.

Ph.\ Barbe, W.P.\ McCormick (2004).  Asymptotic expansions for
infinite weighted convolutions of heavy tail distributions and
applications, {\tt http://www.arxiv.org/abs/math.PR/0412537}, submitted.

N.H.\ Bingham, C.M.\ Goldie, J.L.\ Teugels (1989). {\sl Regular Variation}, 
2nd ed. Cambridge University Press.

Y.S.\ Chow, H.\ Teicher (1978). {\sl Probability Theory, Independence, 
Interchangeability, Martingales}, Springer.

K.L.\ Chung (1974). {\sl A Course in Probability Theory}, 2nd ed., 
Academic Press.

R.A.\ Doney (1980). Moments of Ladder Heights in
Random Walks, {\sl J.\ Appl.\ Prob.}, 17, 248--252.

P.\ Embrechts, C.M.\ Goldie, N.\ Veraverbeke (1979). Subexponentiality
and infinite divisibility, {\sl Z.\ Wahrsch.\ verw.\ Geb.}, 49, 335--347. 

P.\ Embrechts, C.\ Kl\"uppelberg, T.\ Mikosch (1997). {\sl Modelling 
Extremal Events}, Springer.

J.L.\ Geluk (1992). Second order tail behaviour of
a subordinated probability distribution. {\sl Stoch. Proc. Appl.},
40, 325--337.

J.L.\ Geluk (1996). Tails of subordinated laws: The
regularly varying case. {\sl Stoch. Proc. Appl.}, 61, 147--161.

P.\ Embrechts, N.\ Veraverbeke, N. (1982). Estimates for 
the probability of ruin with special emphasis on
the possibility of large claims, {\sl Insurance: Mathematics and
Economics}, 1, 55--72.

W.\ Feller (1971). {\sl An Introduction to
Probability Theory and its Applications}, 2nd ed., Wiley.

S.\ Karlin, H.M.\ Taylor (1975). {\sl A First Course in Stochastic
Processes}, 2nd ed., Academic Press.

V.~Mari\'c (2000). {\sl Regular Variation and Differential Equations},
{\sl Lecture Notes in Mathematics}, 1726, Springer.

E.\ Omey, E.\ Willekens (1986). Second order
behaviour of the tail of a subordinated probability distribution.
{\sl Stoch. Proc. Appl.}, 2, 339--353.

E.\ Omey, E.\ Willekens (1987). Second-order
behaviour of distributions subordinate to a distribution with finite
mean.  {\sl Comm. Statist. Stoch. Models}, 3, 311--342.

N.\ Veraverbeke (1977). Asymptotic behavior of
Wiener-Hopf factors of a random walk. {\sl Stoch. Proc. Appl.}, 5, 27--37.

}

\vskip .5in

\setbox1=\vbox{\halign{#\hfil&\hskip 40pt #\hfill\cr
  Ph.\ Barbe            & W.P.\ McCormick and C.\ Zhang\cr
  90 rue de Vaugirard   & Dept.\ of Statistics \cr
  75006 PARIS           & University of Georgia \cr
  FRANCE                & Athens, GA 30602 \cr
                        & USA \cr
                        & $\{$bill,czhang$\}$@stat.uga.edu \cr}}
\box1

\vfill
\bye